\documentclass[10pt,twocolumn,twoside]{IEEEtran}

\usepackage{bm}   
\usepackage{amsmath}
\usepackage{amssymb}
\usepackage{indentfirst}  
\usepackage{graphicx}  
\usepackage{multirow} 
\usepackage{booktabs} 
\usepackage{cite}  
\usepackage{threeparttable}   
\usepackage{algorithm}
\usepackage{algpseudocode}
\usepackage{arydshln}
\usepackage{color}
\usepackage[american]{babel}
\usepackage{microtype}
\usepackage[numbers,sort&compress]{natbib}  
\usepackage{steinmetz}

\usepackage{cases}
\usepackage{enumerate}
\usepackage{subcaption} 
\usepackage{epstopdf}
\usepackage{url}

\allowdisplaybreaks[1] 

\newtheorem{theorem}{Theorem}
\newtheorem{corollary}{Corollary}
\newtheorem{lemma}{Lemma}
\newtheorem{definition}{Definition}
\newtheorem{remark}{Remark}

\def\mV{\mathcal{V}}
\def\mE{\mathcal{E}}

\def\mG{\mathcal{G}}

\def\mbR{\mathbb{R}}
\def\mbC{\mathbb{C}}

\def\Lg{\bm{L}}

\begin{document}
\title{On Directed Graphs With Real Laplacian Spectra}
	
\author{
		Tianhao~Yu,~
		Shenglu~Wang,~
		Mengqi~Xue,~\IEEEmembership{Member,~IEEE,}
		Yue~Song,~\IEEEmembership{Member,~IEEE,}
		and~David~J.~Hill,~\IEEEmembership{Life Fellow,~IEEE}
		
\thanks{T. Yu, S. Wang, M. Xue and Y. Song are with the Department of Control Science and Engineering, Tongji University, Shanghai 201804, China, also with the State Key Laboratory of Autonomous Intelligent Unmanned Systems, Shanghai 201210, China, and also with the Frontiers Science Center for Intelligent Autonomous Systems, Ministry of Education, Shanghai 200120, China (e-mail: ysong@tongji.edu.cn).}
\thanks{D. J. Hill is with the Department of Electrical and Computer Systems Engineering, Monash University, Clayton, VIC 3800, Australia, and also with the Department of Electrical and Electronic Engineering, University of Hong Kong, Hong Kong (e-mail: davidj.hill@monash.edu).}
}
	
\markboth{}  
{Yu \MakeLowercase{\textit{et al.}}: Digraphs With Real Laplacian Spectra}

\maketitle

\begin{abstract}
	
	It is reported that dynamical systems over digraphs have superior performance in terms of system damping and tolerance to time delays if the underlying graph Laplacian has a purely real spectrum. This paper investigates the topological conditions under which digraphs possess real or complex Laplacian spectra. We derive sufficient conditions for digraphs, which possibly contain self-loops and negative-weighted edges, to have real Laplacian spectra. The established conditions generally imply that a real Laplacian spectrum is linked to the absence of the so-called digon sign-asymmetric interactions and non-strong connectivity in any subgraph of the digraph. Then, two classes of digraphs with complex Laplacian spectra are identified, which imply that the occurrence of directed cycles is a major factor to cause complex Laplacian eigenvalues. Moreover, we extend our analysis to multilayer digraphs, where strategies for preserving real/complex spectra from graph interconnection are proposed. Numerical experiments demonstrate that the obtained results can effectively guide the redesign of digraph topologies for a better performance.
\end{abstract}


%
\IEEEpeerreviewmaketitle

\section{Introduction}

In multi-agent systems (MASs), directed interactions are ubiquitous, shaping complex dynamics that underpin applications ranging from autonomous vehicle platoons to smart grids \cite{zheng2015stability, gan2014exact}. The underlying network topology of an MAS is naturally encoded in the Laplacian matrix, which is a key analysis tool in the study of the system collective behaviors such as consensus, synchronization, stability and performance metrics \cite{olfati2007consensus, qin2020synchronization, pirani2019graph, dezfulian2021performance}. Unlike undirected graphs, which have symmetric Laplacian matrices with purely real spectra, the Laplacian matrices of directed graphs (also known simply as digraphs) may contain complex eigenvalues. However, the presence of complex Laplacian eigenvalues can cause multiple adverse effects on the system performance. Firstly, in the consensus problem of MASs, the damping ratio of the system will decrease due to the nonzero imaginary part of the Laplacian eigenvalues, resulting in slower decay and lasting oscillations \cite{chowdhury2020scalable}. Secondly, the studies on several dynamical network models \cite{bliman2008average, dezfulian2021performance, ma2022delay, moradian2018robustness} reveal that the presence of complex Laplacian eigenvalues makes the consensus or synchronization less robust to the time delays in the agents' interactions, which is inevitable in many real-world MASs. Thirdly, the definition of algebraic connectivity, which is a fundamental concept in graph theory, becomes indefinite with the presence of complex Laplacian eigenvalues. For undirected graphs, the algebraic connectivity is defined as the second smallest Laplacian eigenvalue \cite{fiedler1973algebraic}. In case of complex Laplacian eigenvalues, the concept of the second smallest eigenvalues becomes indefinite. Some works define the algebraic connectivity by the eigenvalue of the digraph Laplacian with the second smallest real part \cite{asadi2016generalized, gao2025effects}, while other works define it by the second smallest eigenvalue of a real symmetric matrix induced from the digraph Laplacian \cite{olfati2004consensus, yu2009second}. All the above issues bring the critical question to our attention: under which condition does or does not a digraph have complex Laplacian eigenvalues?

This question has been ignored over a long period of time as it has been generally assumed that the Laplacian matrix of a digraph contains complex eigenvalues in the analysis of engineering problems \cite{li2009consensus, mei2015distributed, chen2022leader, gao2025distributed}. But it can be observed from numerical experiments that not all digraphs contain complex numbers in their Laplacian spectra. It is of interest to specify those classes of digraphs with real Laplacian spectra since the dynamical systems over these digraphs have better performances in several aspects as aforementioned.

There have been few works investigating the relationship between digraph topologies and spectral complex-valuedness of graph-induced matrices. Ref. \cite{van2018directed} focuses on the addition of a single directed edge to an undirected graph and discovers that the adjacency matrix will have complex eigenvalues if the weight of the directed edge is greater than a threshold determined by the original adjacency matrix. It is revealed in \cite{ochi2022contribution} that if the digraph is obtained by applying a uniformly random edge directization to an undirected graph, then the relative positions of the real part of the adjacency matrix spectrum is conserved in the sense of first-order expansion. When it comes to Laplacian matrices, \cite{ahmadizadeh2017eigenvalues} studies the impact of negative weighted edges on the eigenvalues of digraph Laplacian matrices and provides conditions under which the real part of the non-zero Laplacian eigenvalues remain positive. It is proved in a recent work \cite{cao2025synchronization} that for a special kind of single-root digraphs, the Laplacian eigenvalue with the second smallest real part will be a purely real number.  Nevertheless, it remains open to identify digraphs with entirely real Laplacian spectra in terms of their topological features.

This paper aims to figure out those key topological factors in digraphs which leads to real or complex Laplacian spectra. The contributions are threefold:

\begin{enumerate}
	\item We derive sufficient conditions for digraphs (possibly with self-loops and negative-weighted edges) to have real Laplacian spectra, which require the absence of digon sign-asymmetric interactions (i.e., the bidirectional edges between a pair of nodes have opposite signs) and non-strong connectivity in any subgraph of the digraph.
	\item We identify two specific classes of digraphs that must have complex Laplacian eigenvalues. A common property shared by these two classes of digraphs is the existence of directed cycles.
	\item The obtained results are extended to multilayer digraphs, which are obtained by interconnecting some individual digraphs as the components. We propose some strategies for constructing multilayer digraphs by which the real or complex Laplacian spectrum of certain graph components is transferred to the multilayer one. This result provides a guideline for digraph expansion with the preservation of a real Laplacian spectrum.
\end{enumerate}

The rest of the paper is organized as follows. The notations and preliminaries used in this paper are given in Section \ref{sec:notations}. In Section \ref{sec:real_Laplacian_spectra}, the sufficient conditions for digraphs to have purely real Laplacian spectra are presented. Section \ref{sec:complex_Laplacian_spectra} characterizes digraphs with complex Laplacian spectra. Section \ref{sec:multilayer_digraph} extends the analysis to multilayer digraphs. Section \ref{sec:Numerical_examples} illustrates the obtained results through numerical examples, and Section \ref{sec:conclusion} gives a conclusion of this paper and an outlook on future directions.

\section{Notations and preliminaries}\label{sec:notations}
Let $\mbC$ denote the set of complex numbers and $\mbR$ the set of real numbers. The square root of $-1$ is denoted by $\imath$. The notation $\bm{I}_p\in \mbR^{p \times p}$ denotes an identity matrix, $\bm{1}_{p \times q}\in \mbR^{p \times q}$ denotes a matrix with all entries being one, $\bm{0}_{p \times q}\in \mbR^{p \times q}$ denotes the zero matrix.

A weighted digraph $\mG$ is denoted by the tuple $\mG(\mV,\mE,\bm{A})$ where $\mV$ is the set of nodes with $n = |\mV|$, $\mE\subseteq\mV\times\mV$ is the set of directed edges, $\bm{A}=[a_{ij}]\in\mbR^{|\mV|\times |\mV|}$ is the weighted adjacency matrix.
An edge oriented from node $j$ to node $i$ is denoted by the \textit{ordered pair} $(i,j)\in\mE$.
The adjacency matrix $\bm{A}$ is defined as follows\footnote{Note that we set $a_{ij}\neq 0$ if there is a directed edge from node $j$ to node $i$, implying that node $i$ can receive information from node $j$. This type of definition is commonly adopted in the area of network consensus and synchronization.}: $a_{ij}\neq0$ denotes the weight of edge $(i,j)\in\mE$, $a_{ij}=0$ if $(i,j)\notin\mE$, $a_{ii}=0$ if there is no self-loop at node $i$ and $a_{ii}\neq 0$ denotes the weight of the self-loop at node $i$.
A graph $\mG(\mV,\mE,\bm{A})$ is \textit{unweighted} if $a_{ij}=1$, $\forall (i,j)\in\mE$.
A graph $\mG(\mV,\mE,\bm{A})$ is \textit{loopy} if there exists $a_{ii}\neq 0$, or \textit{loopless} if $a_{ii}=0$, $\forall i\in\mV$.
The Laplacian matrix $\Lg=[L_{ij}]\in\mbR^{n\times n}$ is defined as $L_{ij}=-a_{ij}$, $L_{ii}=\sum_{j=1}^n a_{ij}$.
Note that negative edge weights $a_{ij}<0$ are allowed in this paper. The Laplacian matrix of a graph with negative edge weights is still possible to be positive semi-definite, referring to some recent literature \cite{song2019extension,mukherjee2019robustness} for the details.
Given a digraph $\mG(\mV,\mE,\bm{A})$, $\mG_1(\mV_1,\mE_1,\bm{A}_1)$ is its subgraph where $\mV_1\subseteq\mV$ and $(i,j)\in\mE_1$ if $i,j\in\mV_1$ and $(i,j)\in\mE$.
For simplicity, we use $\bm{L}_{\mV_1,\mV_2}$ to denote the block submatrix of $\bm{L}$ whose rows are indexed by the subset of nodes $\mV_1$ and columns are indexed by the subset of nodes $\mV_2$. Note that $\bm{L}_{\mV_1,\mV_1}$ represents the Laplacian matrix of the subgraph $\mG_1(\mV_1,\mE_1,\bm{A}_1)$ with those incoming connections from $\mV\backslash \mV_1$ to $\mV_1$ being changed to self-loops at the associated nodes in $\mV_1$.

\indent
A \textit{directed path} in a digraph refers to an ordered sequence of nodes $\{i_1,i_2,...,i_p\}$ such that any pair of consecutive nodes in the sequence $(i_k,i_{k+1})$ is an edge of the digraph.
A \textit{directed cycle} in a digraph refers to a directed path linking an ordered sequence nodes of at least three nodes, such that the last node coincide with the first node.
A digraph is \textit{strongly connected} if there exists a directed path from any node to any other node.
A digraph $\mG(\mV,\mE,\bm{A})$ is actually undirected if $\bm{A}=\bm{A}^T$.
The \textit{undirected version} of a digraph $\mG(\mV,\mE,\bm{A})$ is obtained by ignoring the directions of all edges $\mE$.

\section{Digraphs with real Laplacian spectra}\label{sec:real_Laplacian_spectra}
In a digraph $\mG(\mV,\mE,\bm{A})$, two edges linking the same pair of nodes from two orientations (i.e., $(i,j)\in\mE$ and $(j,i)\in\mE$) form a digon.
To better present our results, we first introduce the following definitions about the digon-type interactions.
\begin{definition}\label{def:digon}
	Given a digraph $\mG(\mV,\mE,\bm{A})$, the interaction between node $i$ and node $j$ is said to be
	\begin{itemize}
		\item \textit{unidirectional} if $a_{ij} \neq 0, a_{ji}=0$ or $a_{ij}=0, a_{ji} \neq 0$
		\item \textit{digon symmetric} if $a_{ij}=a_{ji}\neq0$;
		\item \textit{digon asymmetric} if $a_{ij}\neq 0$, $a_{ji}\neq0$ and $a_{ij}\neq a_{ji}$;
		\item \textit{digon sign-asymmetric} if $a_{ij}a_{ji}<0$.
	\end{itemize}
\end{definition}

Definition \ref{def:digon} classifies the possible interactions between a pair of nodes. A digon symmetric interaction actually implies an undirected edge between the two nodes. Digon asymmetric and digon sign-asymmetric respectively refer to the different weights and signs of interaction between two nodes. An unidirectional interaction implies a single directed edge between the two nodes. These concepts play an important role in specifying digraphs with purely real Laplacian spectra.

The following lemma characterizes the spectral property held by two-node digraphs, which is simple but helpful in the subsequent analysis.

\begin{lemma}\label{lm:two_node}
	For any two-node digraph that contains no digon sign-asymmetric interactions and possibly self-loops, the eigenvalues of its Laplacian matrix are all real numbers.
\end{lemma}

\begin{IEEEproof}
	The Laplacian matrix of such a two-node digraph can be expressed as:
	$$ \Lg = \left[\begin{matrix}
		a_{11} & a_{12} \\
		a_{21} & a_{22}
	\end{matrix}\right].$$
	Its characteristic equation is quadratic
	$\lambda^2 - (a_{11}+a_{22})\lambda + a_{11}a_{22} - a_{12}a_{21} = 0$
	with a positive discriminant $\Delta = (a_{11}-a_{22})^2 + 4a_{12}a_{21} \geq 0$,
	which completes the proof.
\end{IEEEproof}

\begin{lemma}\label{lm:real_numbers}
	The eigenvalues of the Laplacian matrix of a possibly loopy digraph are all real numbers if any of its subgraphs is not strongly connected, or actually undirected.
\end{lemma}

\begin{IEEEproof}
	According to \cite{bullo2018lectures}, if a digraph $\mG(\mV,\mE,\bm{A})$ is not strongly connected, then its Laplacian matrix is reducible and hence can be transformed into an upper block triangular matrix by renumbering the nodes in the $\mG$, as shown below:
	\begin{equation}\label{eq:block_L}
		\Lg = \left[\begin{matrix}
			\Lg_{\mV_1, \mV_1} & \Lg_{\mV_1, \mV \setminus \mV_1} \\
			\bm{0} & \Lg_{\mV \setminus \mV_1, \mV \setminus \mV_1}
		\end{matrix}\right]
	\end{equation}
	where $\Lg_{\mV_1, \mV_1}$ and $\Lg_{\mV \setminus \mV_1, \mV \setminus \mV_1}$, as defined in Section \ref{sec:notations}, represent the Laplacian matrices of the subgraphs $\mG_1$ and $\mG \setminus \mG_1$ with incoming connections changed to self-loops. Subgraphs $\mG_1$ and $\mG\setminus\mG_1$ are only connected via those unidirectional interactions indexed by those nonzero entries in $\Lg_{\mV_1, \mV \setminus \mV_1}$. According to the block structure in (\ref{eq:block_L}), we have $\textup{spec}(\Lg)=\textup{spec}(\Lg_{\mV_1, \mV_1}) \cup \textup{spec}(\Lg_{\mV \setminus \mV_1, \mV \setminus \mV_1})$, where $\textup{spec}()$ denotes the matrix spectrum.
	
	According to the condition in Lemma \ref{lm:real_numbers}, $\mG_1$ is undirected or not strongly connected. In the former case, $\Lg_{\mV_1, \mV_1}$ is real symmetric with a real Laplacian spectrum. In the latter case, subgraph $\mG_1$ is still not strongly connected, and thus we can iteratively find the matrix partition similar to (\ref{eq:block_L}) and finally $\Lg_{\mV_1, \mV_1}$ takes the following partitioned form:
	\begin{multline}\label{eq:multi_block_L}
		\Lg_{\mV_1, \mV_1} = \\
		\left[\begin{array}{c:c}
			\begin{array}{c:c}
				\begin{array}{c:c}
					\Lg_{\mV_l, \mV_l} & \cdots \\
					\hdashline
					\vdots & \ddots
				\end{array} & \Lg_{\mV_3, \mV_2 \setminus \mV_3} \\
				\hdashline
				\bm{0} & \Lg_{\mV_2 \setminus \mV_3, \mV_2 \setminus \mV_3}
			\end{array} & \Lg_{\mV_2, \mV_1 \setminus \mV_2}\\
			\hdashline
			\bm{0} & \Lg_{\mV_1 \setminus \mV_2, \mV_1 \setminus \mV_2}
		\end{array}\right]
	\end{multline}
	where the each principal submatrix from $\Lg_{\mV_1 \setminus \mV_2, \mV_1 \setminus \mV_2}$ to $\Lg_{\mV_l, \mV_l}$ is the Laplacian matrix of a possibly loopy undirected or single-node graph. Note that the above decomposition applies to $\Lg_{\mV \setminus \mV_1, \mV \setminus \mV_1}$ as well, and hence it follows that
	\begin{align*}
		\textup{spec}(\Lg)=&\textup{spec}(\Lg_{\mV_l, \mV_l}) \cup \cdots \cup \textup{spec}(\Lg_{\mV_2 \setminus \mV_3, \mV_2 \setminus \mV_3}) \\
		\cup &\textup{spec}(\Lg_{\mV_1 \setminus \mV_2, \mV_1 \setminus \mV_2}) \cup \textup{spec}(\Lg_{\mV \setminus \mV_1, \mV \setminus \mV_1})
	\end{align*}
	is a set of real numbers.
\end{IEEEproof}

Note that any subgraph being not strongly connected, as required in Lemma \ref{lm:real_numbers}, implies that the interaction between each pair of nodes can only be unidirectional, which is rather restrictive. Combining Lemma \ref{lm:two_node} and Lemma \ref{lm:real_numbers}, we are ready to specify a more general class of digraphs the Laplacian matrices of which have purely real spectra.

\begin{theorem}\label{th:real_numbers}
	The eigenvalues of the Laplacian matrix of a possibly loopy digraph are all real numbers if the digraph satisfies both the following conditions:
	\begin{enumerate}
		\item it does not contain digon sign-asymmetric interactions;
		\item any of its subgraphs with three or more nodes is not strongly connected or actually undirected.
	\end{enumerate}
\end{theorem}

\begin{IEEEproof}
	Under the condition in Theorem \ref{th:real_numbers}, the iterative matrix partition in (\ref{eq:multi_block_L}) in the proof of Lemma \ref{lm:real_numbers} still applies to $\Lg$, which will end when each principal submatrix is the Laplacian matrix of one of the following three types of possibly loopy graphs:
	\begin{enumerate}
		\item an undirected graph
		\item a single-node graph
		\item a two-node digraph without digon sign-asymmetric interactions
	\end{enumerate}
	all of which have been solved in Lemma \ref{lm:two_node} and Lemma \ref{lm:real_numbers}.
\end{IEEEproof}

\indent
Theorem \ref{th:real_numbers} implies that realness of the spectrum of a Laplacian matrix of a digraph is closely linked to the property of non-strong connectivity. The proof of Theorem \ref{th:real_numbers} is essentially an iterative decomposition of the Laplacian matrix of the digraph into the Laplacian matrices of some single-node, two-node or undirected subgraphs, which leads to a real Laplacian spectrum. Recall that \cite{cao2025synchronization} studies the complex-valuedness of the Laplacian eigenvalue with the second smallest real part and proves that this particular eigenvalue is a real number if the digraph is positive weighted, loopless and contains a single root and a spanning tree. By comparison, Theorem \ref{th:real_numbers} gives a topological description of digraphs with entirely real spectra, while taking self-loops and negative-weighted edges into consideration. This result reveals that it is not difficult to construct a digraph with a real Laplacian spectrum.

\indent
Theorem 1 also leads to the following corollary which specifies a more illustrative type of digraphs with purely real Laplacian spectra.

\begin{definition}
	A digraph is called \textit{tree-type digraph} if its undirected version is a tree.
\end{definition}

Note that the tree-type digraph has a more relaxed description than the concept of directed tree \cite{bullo2018lectures}, which is an acyclic digraph containing a unique root node such that any other node of the digraph can be reached by a unique directed path starting at the root.

\begin{corollary}\label{coro:tree}
	For a digraph without digon asymmetric interactions, if it is a tree-type digraph, then the eigenvalues of the Laplacian matrix of the digraph are all real numbers.
\end{corollary}

\begin{IEEEproof}
	A tree-type digraph without digon asymmetric interactions means this digraph only contains unidirectional or digon symmetric interactions, so it falls into the description of Theorem \ref{th:real_numbers}.
\end{IEEEproof}

\indent
The tree-type topology specified in Corollary \ref{coro:tree} has demonstrated significant practicality, especially in MASs. For instance, it is well known in vehicle platoons since it is easy to construct and beneficial to the performance of vehicle platoons \cite{pirani2022impact, zheng2015stability}. Moreover, it is widely used in smart grids, and many studies are based on it to solve control or optimization problems in smart grids \cite{gan2014exact, li2017coordinated}.

\section{Digraphs with complex Laplacian spectra}\label{sec:complex_Laplacian_spectra}

This section will specify some classes of digraphs that surely have complex Laplacian spectra. Let us begin with the unweighted directed cycle graph, which is common and representative in practice. Moreover, since the Laplacian matrix of an unweighted directed cycle graph is well-structured, an analytical expression for its eigenvalues is available.

\begin{theorem}\label{th:cycle}
	The Laplacian matrix of an unweighted directed cycle graph must contain complex eigenvalues.
\end{theorem}

\begin{IEEEproof}
	The characteristic equation of $\Lg$ of an unweighted directed cycle graph is:
	\begin{align*}
		\left|\Lg - \lambda \bm{I} \right| &= \left|\begin{matrix}
			1-\lambda & -1 & \cdots & 0\\
			0 & 1-\lambda & \cdots & 0\\
			\vdots & \vdots &  \ddots & \vdots\\
			-1 & 0 & \cdots & 1-\lambda
		\end{matrix}\right| \\
		&= (1 - \lambda)^n - 1 = 0
	\end{align*}
	in which the determinant is expanded along the first column. Then we can derive the expression for the eigenvalues:
	\begin{equation*} \label{eq:cycle}
		\lambda=1-e^{\imath \frac{2 k \pi}{n}}, \quad k=0,1,2, \ldots, n-1
	\end{equation*}
	which must contain complex eigenvalues.
\end{IEEEproof}

\begin{figure}[htbp]
	\centering
	\begin{subfigure}{0.23\textwidth}
		\centering
		\includegraphics[width=0.6\textwidth]{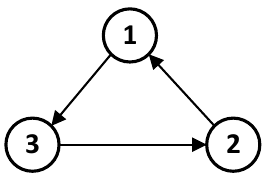}
		\caption{Unweighted.}
		\label{fig:3node_cycle_a}
	\end{subfigure}
	\hfill
	\begin{subfigure}{0.23\textwidth}
		\centering
		\includegraphics[width=0.6\textwidth]{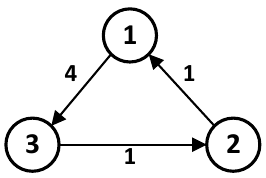}
		\caption{With edge weights $1$, $1$, $4$.}
		\label{fig:3node_cycle_b}
	\end{subfigure}
	\caption{Examples of three-nodes directed cycle graphs.}
	\label{fig:3node_cycle}
\end{figure}

\begin{remark}
	Note that a weighted directed cycle graph may have a purely real Laplacian spectrum. For instance, for a three-node unweighted directed cycle (shown in Fig. \ref{fig:3node_cycle}(\subref{fig:3node_cycle_a})), the eigenvalues of its Laplacian matrix are $0$, $1.5 \pm 0.87 \imath$ (to two decimal places), but for a three-node directed cycle with weights $1$, $1$, $4$ (shown in Fig. \ref{fig:3node_cycle}(\subref{fig:3node_cycle_b})), the eigenvalues of its Laplacian matrix are $0$, $3$, $3$.
\end{remark}

The next definition can help us describe a kind of more complicated loopless digraph containing a directed cycle, and prove by Theorem \ref{th:complete_cycle} below that its Laplacian matrix must contain complex eigenvalues.

\begin{definition}\label{def:directed_cycle-embedded_complete_graph}
	A digraph is called an \textit{unweighted directed cycle-embedded complete (UDC-EC) graph} if it is constructed in the following way. Given a $n$-node unweighted undirected complete graph, then select some pairs of nodes and change the originally digon symmetric interactions between these pairs of nodes to unidirectional interactions such that all these unidirectional interactions form a directed cycle.
\end{definition}

\begin{figure}[htbp]
	\centering
	\includegraphics[width=0.21\textwidth]{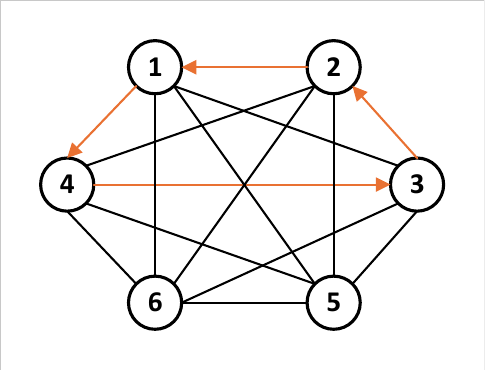}
	\caption{An example of a six-node UDC-EC graph.}
	\label{fig:directed_cycle-embedded_complete_graph}
\end{figure}

An example of a six-node UDC-EC graph described by Definition \ref{def:directed_cycle-embedded_complete_graph} is shown in Fig. \ref{fig:directed_cycle-embedded_complete_graph}, which has a total of four unidirectional interactions forming a directed cycle. For simplicity, in Fig. \ref{fig:directed_cycle-embedded_complete_graph}, the undirected interactions are represented by edges without arrows, and the unidirectional interactions are marked in orange.

\begin{theorem}\label{th:complete_cycle}
	The Laplacian matrix of a UDC-EC graph must contain complex eigenvalues.
\end{theorem}

\begin{IEEEproof}
	Assume the UDC-EC graph contains $n$ nodes in total and the directed cycle involves $m$ nodes. Then the Laplacian matrix can be expressed as:
	\begin{equation*}\label{eq:complete_cycle_block_L}
		\Lg = \left[\begin{matrix}
			\Lg_{m} & -\bm{1} \\
			-\bm{1} & \Lg_{n-m}
		\end{matrix}\right]
	\end{equation*}
	where $\Lg_{m} \in \mbR^{m \times m}$ and $\Lg_{n-m} \in \mbR^{(n-m) \times (n-m)}$ are
	\begin{equation*}
		\Lg_{m} = \left[\begin{matrix}
			n-2 & -1 & -1 & \cdots & -1 & 0\\
			0 & n-2 & -1 & \cdots & -1 & -1\\
			-1 & 0 & n-2 & \cdots & -1 & -1\\
			\vdots & \vdots & \vdots & \ddots & \vdots & \vdots\\
			-1 & -1 & -1 & \cdots & n-2 & -1\\
			-1 & -1 & -1 & \cdots & 0 & n-2\\
		\end{matrix}\right],
	\end{equation*}
	\begin{equation*}
		\Lg_{n-m} = \left[\begin{matrix}
			n-1 & -1 & \cdots & -1\\
			-1 & n-1 & \cdots & -1\\
			\vdots & \vdots & \ddots & \vdots\\
			-1 & -1 & \cdots & n-1\\
		\end{matrix}\right].
	\end{equation*}
	
	We still derive the expression for the eigenvalues by solving the characteristic equation $\left| \Lg - \lambda \bm{I} \right| = 0$. The following determinant transformations are applied in turn to simplify the characteristic polynomial:
	\begin{enumerate}
		\item Add the second to $n^{th}$ columns to the first column in sequence, making all entries of the first column $-\lambda$.
		\item Factor out $-\lambda$ from the first column of the determinant, resulting in all entries of that column being 1.
		\item Add the first column to each of the rest of the columns.
	\end{enumerate}
	After that, the characteristic equation can be expressed in the following lower block triangular form:
	\begin{equation}\label{eq:complete_cycle_block_D}
		-\lambda \left|\begin{matrix}
			\bm{D}_m & \bm{0} \\
			\bm{d} & \bm{D}_{n-m}
		\end{matrix}\right|= -\lambda \bm{D}_m \bm{D}_{n-m}=0
	\end{equation}
	where
	\begin{equation}\label{eq:Dm}
		\bm{D}_m = \left|\begin{matrix}
			1 & 0 & 0 & \cdots & 0 & 1\\
			1 & \delta-1 & 0 & \cdots & 0 & 0\\
			1 & 1 & \delta-1 & \cdots & 0 & 0\\
			\vdots & \vdots & \vdots & \ddots & \vdots & \vdots\\
			1 & 0 & 0 & \cdots & \delta-1 & 0\\
			1 & 0 & 0 & \cdots & 1 & \delta-1\\
		\end{matrix}\right|, \tag{\ref{eq:complete_cycle_block_D}a}
	\end{equation}
	\begin{equation}\label{eq:Dn-m}
		\bm{D}_{n-m} = \left|\begin{matrix}
			\delta & 0 & \cdots & 0\\
			0 & \delta & \cdots & 0\\
			\vdots & \vdots & \ddots & \vdots\\
			0 & 0 & \cdots & \delta\\
		\end{matrix}\right|, \tag{\ref{eq:complete_cycle_block_D}b}
	\end{equation}
	$\delta = n - \lambda$ and $\bm{d} = \left|\bm{1}, \bm{0}, \cdots, \bm{0}\right| \in \mbR^{(n-m) \times m}$.
	
	It can be easily observed from (\ref{eq:complete_cycle_block_D}) and (\ref{eq:Dn-m}) that $\Lg$ has a simple zero eigenvalue and $n-m$ eigenvalues of value $n$. For (\ref{eq:Dm}), we expand it along the first row:
	\begin{equation}\label{eq:D_m}
		\bm{D}_m = (\delta - 1)^{m-1} + (-1)^{m+1} \bm{D}^\prime_{m-1}  = 0
	\end{equation}
	where
	$$\bm{D}^\prime_{m-1} = \left|\begin{matrix}
		1 & \delta - 1 & 0 & \cdots & 0 & 0\\
		1 & 1 & \delta - 1 & \cdots & 0 & 0\\
		1 & 0 & 1 & \cdots & 0 & 0\\
		\vdots & \vdots & \vdots & \ddots & \vdots & \vdots\\
		1 & 0 & 0 & \cdots & 1 & \delta - 1 \\
		1 & 0 & 0 & \cdots & 0 & 1
	\end{matrix}\right|.$$
	Then, expand $\bm{D}^\prime_{m-1}$ along the last row:
	$$\bm{D}^\prime_{m-1} = (-1)^{m}(\delta - 1)^{m-2} + \bm{D}^\prime_{m-2} .$$
	Since $\bm{D}^\prime_{m-2}$ is structurally similar to $\bm{D}^\prime_{m-1}$, by expanding $\bm{D}^\prime_{m-k}\quad(k = 2, 3, \ldots)$, we obtain:
	\begin{equation}\label{eq:D'm-1}
		\bm{D}^\prime_{m-1} = \frac{1-(-1)^{m-1}(\delta - 1)^{m-1}}{\delta}.
	\end{equation}
	Substituting (\ref{eq:D'm-1}) into (\ref{eq:D_m}) gives:
	$$\bm{D}_m = \frac{(-1)^{m+1} + (\delta - 1)^m}{\delta} = 0.$$
	So the expression of $\lambda$ is $\lambda = n-1+e^{\imath \frac{2 k \pi}{m}}$, where $k=1,2, \ldots, m-1$. In summary, the expression for the Laplacian eigenvalues of the digraph is:
	\begin{equation*}\label{eq:complete_cycle}
		\lambda=
		\begin{cases}
			0, \quad k=0 \\
			n-1+e^{\imath \frac{2 k \pi}{m}}, \quad k=1,2, \ldots, m-1 \\
			n, \quad k= m, \ldots, n-1 \\
		\end{cases}
	\end{equation*}
	which must contain complex eigenvalues.
\end{IEEEproof}


\section{Results on multilayer digraphs}\label{sec:multilayer_digraph}
This section move on to the Laplacian spectrum of a multilayer digraph obtained by interconnecting several subgraphs.
	Multilayer digraphs are a common model to capture the multilayer structure of a variety of real-world complex systems such as social networks, transportation networks and cyber-physical networks \cite{della2020symmetries, de2023more}.
	
	\indent
	Let us first investigate how to interconnect two digraphs to construct a multilayer digraph with a real Laplacian spectrum.
	Consider a multilayer digraph $\mG_3(\mV_1\cup\mV_2,\mE_1\cup\mE_2\cup\mE_{12}\cup\mE_{21})$ composed of two subgraphs $\mG_1(\mV_1,\mE_1,\bm{A}_1)$,  $\mG_2(\mV_2,\mE_2,\bm{A}_2)$, additional directed edges $\mE_{12} \subseteq \mV_1\times \mV_2$ orienting from some nodes in $\mV_1$ to some nodes in $\mV_2$, and additional directed edges $\mE_{21} \subseteq \mV_2\times \mV_1$ orienting from some nodes in $\mV_2$ to some nodes in $\mV_1$. With these notations, we are ready to present the following result.

\begin{corollary}\label{coro:com_real}
	For a multilayer digraph $\mG_3(\mV_1\cup\mV_2,\mE_1\cup\mE_2\cup\mE_{12}\cup\mE_{21})$ composed of $\mG_1(\mV_1,\mE_1,\bm{A}_1)$ and $\mG_2(\mV_2,\mE_2,\bm{A}_2)$, if $\mG_1$ satisfies the conditions in Theorem 1, $\mG_2$ is a digraph with a real Laplacian spectrum, and $\mE_{12} = \phi$, the weights of the edges in $\mE_{21}$ are arbitrary, then the eigenvalues of the Laplacian matrix of multilayer digraph $\mG_3$ are all real numbers.
\end{corollary}

\begin{IEEEproof}	
	The condition $\mE_{12} = \phi$ implies that the directed edges between $\mG_1$ and $\mG_2$ only orient from $\mV_2$ to $\mV_1$, which makes the multilayer digraph not strongly connected. Thus, the Laplacian matrix of $\mG_3$ can be expressed as:
	\begin{equation}\label{eq:com_block}
		\Lg = \left[\begin{matrix}
			\Lg_{\mV_1, \mV_1} & \Lg_{\mV_1, \mV_2} \\
			\bm{0} & \Lg_{\mV_2, \mV_2}
		\end{matrix}\right]
	\end{equation}
	where $\Lg_{\mV_1, \mV_1}$ represents the Laplacian matrix of $\mG_1$ after changing incoming connections to self-loops, $\Lg_{\mV_2, \mV_2}$ denotes the Laplacian matrix of $\mG_2$. It follows from (\ref{eq:com_block}) that, $\textup{spec}(\Lg)=\textup{spec}(\Lg_{\mV_1, \mV_1}) \cup \textup{spec}(\Lg_{\mV_2, \mV_2})$. Theorem \ref{th:real_numbers} guarantees that the eigenvalues of $\Lg_{\mV_1, \mV_1}$ are all real numbers, and the eigenvalues of $\Lg_{\mV_2, \mV_2}$ are all real numbers as well, so the eigenvalues of $\Lg$ are also all real numbers.
\end{IEEEproof}

\begin{figure}[htbp]
	\centering
	\includegraphics[width=0.31\textwidth]{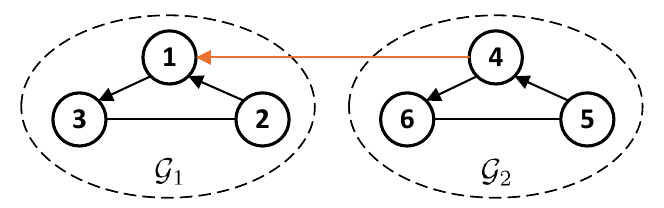}
	\caption{An unweighted multilayer digraph with a complex spectrum while the spectra of $\mG_1$ and $\mG_2$ are both real.}
	\label{fig:multilayer_both_real}
\end{figure}

\begin{remark}
	Note that Corollary \ref{coro:com_real} may not be true if $\mG_1$ only has a real Laplacian spectrum and does not satisfy the conditions in Theorem \ref{th:real_numbers}.
	This is because $\Lg_{\mV_1, \mV_1}$ in (\ref{eq:com_block}) represents the Laplacian matrix of $\mG_1$ after changing the incoming connections from $\mG_2$ into self-loops, by which complex eigenvalues could occur.
	For example, the unweighted multilayer digraph in Fig. \ref{fig:multilayer_both_real} is obtained by interconnecting $\mG_1$ and $\mG_2$.
	The Laplacian spectra of $\mG_1$ and $\mG_2$ are identical and real ($0$, $2$, $2$), but the multilayer graph has a complex Laplacian spectrum ($0.16$, $2.42 \pm 0.61\imath$, $0$, $2$, $2$).
\end{remark}

The next corollary shows how a multilayer graph ``inherits" complex eigenvalues from its subgraph.

\begin{corollary}\label{coro:com_complex}
	For a multilayer digraph $\mG_3(\mV_1\cup\mV_2,\mE_1\cup\mE_2\cup\mE_{12}\cup\mE_{21})$ composed of $\mG_1(\mV_1,\mE_1,\bm{A}_1)$ and $\mG_2(\mV_2,\mE_2,\bm{A}_2)$, if $\mG_1$ is an arbitrary digraph and $\mG_2$ is a digraph with a complex spectrum, and $\mE_{12} = \phi$, the weights of the edges in $\mE_{21}$ are arbitrary, then the Laplacian matrix of $\mG_3$ must contain complex eigenvalues.
\end{corollary}

\begin{IEEEproof}
	Similar to the proof of Corollary \ref{coro:com_real}, the Laplacian matrix of multilayer digraph $\mG_3$ here also takes a block form as in the form of (\ref{eq:com_block}), we have $\textup{spec}(\Lg)=\textup{spec}(\Lg_{\mV_1, \mV_1}) \cup \textup{spec}(\Lg_{\mV_2, \mV_2})$ as well. Thus the complex eigenvalues in $\Lg_{\mV_2, \mV_2}$ lead to the existence of complex eigenvalues in $\Lg$.
\end{IEEEproof}

\begin{figure}[htbp]
	\centering
	\begin{subfigure}{0.23\textwidth}
		\centering
		\includegraphics[width=0.9\textwidth]{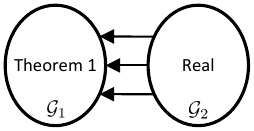}
		\caption{Multilayer digraph with real Laplacian spectrum.}
		\label{fig:illustration_a}
	\end{subfigure}
	\hfill
	\begin{subfigure}{0.23\textwidth}
		\centering
		\includegraphics[width=0.9\textwidth]{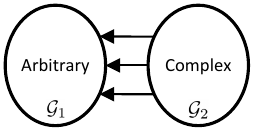}
		\caption{Multilayer digraph with complex Laplacian spectrum.}
		\label{fig:illustration_b}
	\end{subfigure}
	\caption{Illustration of different types of multilayer digraphs.}
	\label{fig:multilayer_illustration}
\end{figure}

Fig. \ref{fig:multilayer_illustration} illustrates the two types of multilayer digraphs corresponding to the description in Corollary \ref{coro:com_real} and Corollary \ref{coro:com_complex}, respectively.
In Fig. \ref{fig:multilayer_illustration}(\subref{fig:illustration_a}), $\mG_1$ is a digraph described in Theorem \ref{th:real_numbers} (marked with ``Theorem 1''), and $\mG_2$ is a digraph with a real Laplacian spectrum (marked with ``Real'').
In Fig. \ref{fig:multilayer_illustration}(\subref{fig:illustration_b}), $\mG_1$ is an arbitrary digraph (marked with ``Arbitrary''), $\mG_2$ is a digraph with a complex Laplacian spectrum (marked with ``Complex'').

The next result illustrates how directed cycle topology affects the spectra of multilayer digraphs.
\begin{definition}\label{def}
	Given an arbitrary graph $\mG(\mV,\mE,\bm{A})$ ($|\mV|=n$), define the multilayer digraph $\mG_m$ as the directed cyclical interconnection of $m$ duplicates of $\mG$ ($m \geq 3$), or \textit{$m$-DCID of $\mG$} in acronym, by the following Laplacian matrix
	\begin{equation}\label{Lap-LWDC}
		\small
		\begin{split}
			\Lg_{\mG_m}=\left[\begin{array}{cccc}
				\Lg_{\mG}+\bm{I}_n & -\bm{I}_n & \cdots & \bm{0} \\
				\bm{0} & \Lg_{\mG}+\bm{I}_n & \cdots & \bm{0} \\
				\vdots & \vdots & \ddots & \vdots \\
				-\bm{I}_n & \bm{0} & \cdots & \Lg_{\mG}+\bm{I}_n
			\end{array}\right]\in\mbR^{(mn\times mn)}
		\end{split}
	\end{equation}
	where $\Lg_{\mG}$ denotes the Laplacian matrix of $\mG$.
\end{definition}

\begin{figure}
	\centering
	\includegraphics[width=0.24\textwidth]{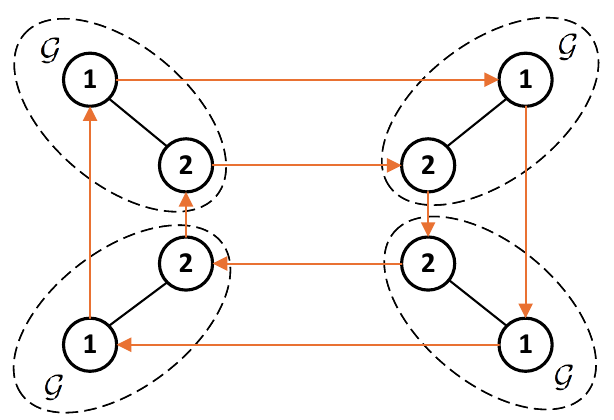}
	\caption{An example of the $4$-DCID of a two-node complete graph.}
	\label{fig:illustration_c}
\end{figure}

\indent
Note that \eqref{Lap-LWDC} is also valid for $m=2$, which refers to an undirected interconnection between two copies of $\mG$. Here we require $m\geq 3$ to construct a proper cycle.
An example of the 4-DCID of a two-node complete graph is shown in Fig.~\ref{fig:illustration_c}.
The following corollary characterizes the Laplacian spectral property of this class of multilayer digraphs.

\begin{corollary}\label{coro:com_real2complex}
	The Laplacian matrix of the $m$-DCID of any graph $\mG$ must contain complex eigenvalues for $m\geq 3$.
\end{corollary}

\begin{IEEEproof}
	The characteristic equation of $\Lg_{\mG_m}$ takes the form
	\begin{equation*}
		\begin{split}
			&|\Lg_{\mG_m} - \lambda \bm{I}|=\\
			&\left|\begin{array}{cccc}
				\Lg_{\mG}+\bm{I}-\lambda\bm{I} & -\bm{I} & \cdots & \bm{0} \\
				\bm{0} & \Lg_{\mG}+\bm{I}-\lambda\bm{I} & \cdots & \bm{0} \\
				\vdots & \vdots & \ddots & \vdots \\
				-\bm{I} & \bm{0} & \cdots & \Lg_{\mG}+\bm{I}-\lambda\bm{I}
			\end{array}\right|=0.
		\end{split}
	\end{equation*}
	Note that the matrix product of every pair of the blocks of $\Lg_{\mG_m} - \lambda \bm{I}$ is commutative. According to Theorem 1 in \cite{silvester2000determinants}, $|\Lg_{\mG_m} - \lambda \bm{I}|$ can be expressed in the form of the Laplace expansion by blocks. Thus, similar to the proof of Theorem \ref{th:cycle}, we have
	\begin{align*}
		|\Lg_{\mG_m} - \lambda \bm{I}|&=|(\Lg_{\mG}+\bm{I}-\lambda\bm{I})^m - \bm{I}|\\
		&=\prod_{i=1}^n\left[\left(\mu_i+1-\lambda\right)^m-1\right]=0
	\end{align*}
	where $\mu_i, i = 1, \dots, n$ are eigenvalues of $\Lg_{\mG}$.
	It leads to the following expression of the eigenvalues of $\Lg_{\mG_m}$
	\begin{equation*}
		\lambda = \mu_i + 1 - e^{\imath \frac{2 k \pi}{m}}, i=1,2,\dots n,~k = 0, 1, \dots, m-1
	\end{equation*}
	which must contain complex values.
\end{IEEEproof}

\indent
Corollary \ref{coro:com_real2complex} is an extended version of Theorem \ref{th:cycle}, which further confirms the presence of directed cycles as a strong indicator of complex Laplacian spectra.
As long as $m$ duplicates of a graph $\mG$ are interconnected in such a way that the backbone structure is a directed cycle, then the resulting multilayer digraph has a complex Laplacian spectrum, regardless of the Laplacian properties of $\mG$.

\section{Numerical examples}\label{sec:Numerical_examples}

In this section, we illustrate the obtained results with several typical examples of digraphs.

\begin{figure}[htbp]
	\centering
	\includegraphics[width=0.3\textwidth]{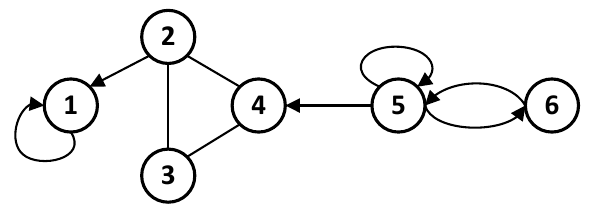}
	\caption{An example of a digraph with real Laplacian spectrum.}
	\label{fig:real_numbers}
\end{figure}

Fig. \ref{fig:real_numbers} shows a six-node digraph of the type described in Theorem \ref{th:real_numbers}. For simplicity, the digon symmetric interactions are represented by undirected edges (without arrows), digon asymmetric but not digon sign-asymmetric interaction is represented by two non-overlapping directed edges, and the node numbering directly derives the upper block triangular matrix form of the Laplacian matrix, as shown in (\ref{eq:example_real}). Suppose $a_{11} = 2$,  $a_{12} = 1$, $a_{23} = a_{32} = 4$, $a_{24} = a_{42} = 1.5$, $a_{34} = a_{43} = 3$, $a_{45} = -7$, $a_{55} = 3.6$, $a_{56} = 1.4$, $a_{65} = 2.1$, so its Laplacian matrix can be expressed as:

\begin{equation}\label{eq:example_real}
	\Lg = \left[\begin{array}{c:ccccc}
		\Lg_{\mV \setminus \mV_1, \mV \setminus \mV_1} & \Lg_{\mV \setminus \mV_1, \mV_1}\\
		\hdashline
		\bm{O} &
		\begin{array}{c:c}
			\Lg_{\mV_1 \setminus \mV_2, \mV_1 \setminus \mV_2}  & \Lg_{\mV_1 \setminus \mV_2, \mV_2} \\
			\hdashline
			\bm{O} & \Lg_{\mV_2, \mV_2}
		\end{array}
	\end{array}\right]
\end{equation}
where $\Lg_{\mV \setminus \mV_1, \mV \setminus \mV_1} = 3,$
$$\Lg_{\mV_1 \setminus \mV_2, \mV_1 \setminus \mV_2} = \left[\begin{array}{ccc}
	5.5 & -4 & -1.5\\
	-4 & 7 & -3\\
	-1.5 & -3 & -2.5\\
\end{array}\right],$$
$$\Lg_{\mV_2, \mV_2} = \left[\begin{array}{cc}
	5 & -1.4\\
	-2.1 & 2.1\\
\end{array}\right],$$
$\Lg_{\mV \setminus \mV_1, \mV_1}$ and $\Lg_{\mV_1 \setminus \mV_2, \mV_2}$ do not affect the eigenvalues of $\Lg$. So the eigenvalues of $\Lg$ are (to two decimal places): $3$, $10.47$, $3.65$, $-4.12$, $5.80$, $1.3$, which do not contain complex numbers.

\begin{figure}[htbp]
	\centering
	\begin{subfigure}{0.23\textwidth}
		\centering
		\includegraphics[width=0.75\textwidth]{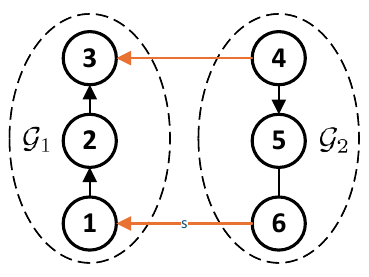}
		\caption{Example of a multilayer digraph with real spectrum.}
		\label{fig:multilayer_real}
	\end{subfigure}
	\hfill
	\begin{subfigure}{0.23\textwidth}
		\centering
		\includegraphics[width=0.75\textwidth]{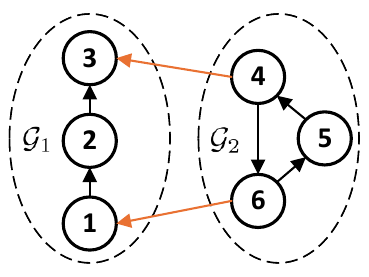}
		\caption{Example of a multilayer digraph with complex spectrum.}
		\label{fig:multilayer_complex}
	\end{subfigure}
	\caption{Examples of multilayer digraphs.}
	\label{fig:multilayer_examples}
\end{figure}

Now we move on to the multilayer graphs. As shown in Fig. \ref{fig:multilayer_examples}(\subref{fig:multilayer_real}), $\mG_1$ satisfies the description of Theorem \ref{th:real_numbers}, $\mG_2$ is a digraph with a real Laplacian spectrum, and the directed edges between them satisfy $\mE_{12} = \phi$, $\mE_{21} \neq \phi$, marked in orange. Let us arbitrarily assign weights to the edges in $\mE_{21}$, e.g.,  $a_{21} = 3$, $a_{32} = 1.8$, $a_{16} = 2$, $a_{34} = 5.3$, $a_{54} = -3.2$, $a_{56} = a_{65} = 1.2$, and direct calculations show that the eigenvalues of Laplacian matrix are $-2.4$, $0$, $1.6$, $2$, $3$, $7.1$, which are all real numbers.

In Fig. \ref{fig:multilayer_examples}(\subref{fig:multilayer_complex}), $\mG_2$ is a three-node unweighted directed cycle, while $\mG_1$ and the edges between $\mG_1$ and $\mG_2$ remain the same as in Fig. \ref{fig:multilayer_examples}(\subref{fig:multilayer_real}). Consequently, the eigenvalues of the Laplacian matrix of this multilayer digraph are $0$, $1.5 \pm \frac{\sqrt{3}}{2} \imath$, $2$, $3$, $7.1$, which contain complex values.


\begin{figure}[htbp]
	\centering
	\begin{subfigure}{0.23\textwidth}
		\centering
		\includegraphics[width=0.6\textwidth]{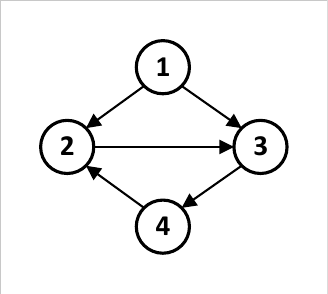}
		\caption{Digraph with a complex Laplacian spectrum.}
		\label{fig:mas1}
	\end{subfigure}
	\hfill
	\begin{subfigure}{0.23\textwidth}
		\centering
		\includegraphics[width=0.6\textwidth]{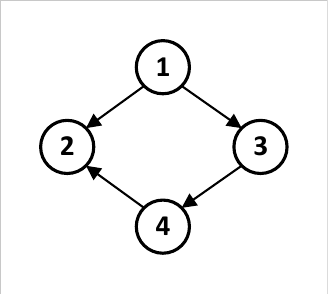}
		\caption{Digraph with a real Laplacian spectrum.}
		\label{fig:mas2}
	\end{subfigure}
	\caption{Examples of MASs with four agents.}
	\label{fig:mas}
\end{figure}

\begin{figure}[htbp]
	\centering
	\begin{subfigure}{0.22\textwidth}
		\centering
		\includegraphics[width=0.99\textwidth]{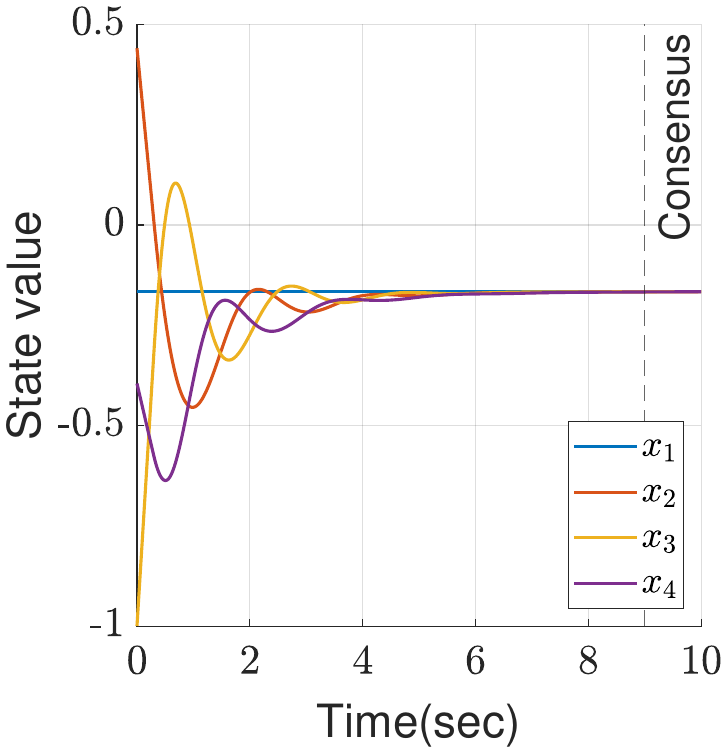}
		\caption{Consensus result of MAS in Fig. \ref{fig:mas}(\subref{fig:mas1}) with $\tau = 0.3s$.}
		\label{fig:consensus_result1}
	\end{subfigure}
	\hfill
	\begin{subfigure}{0.22\textwidth}
		\centering
		\includegraphics[width=0.99\textwidth]{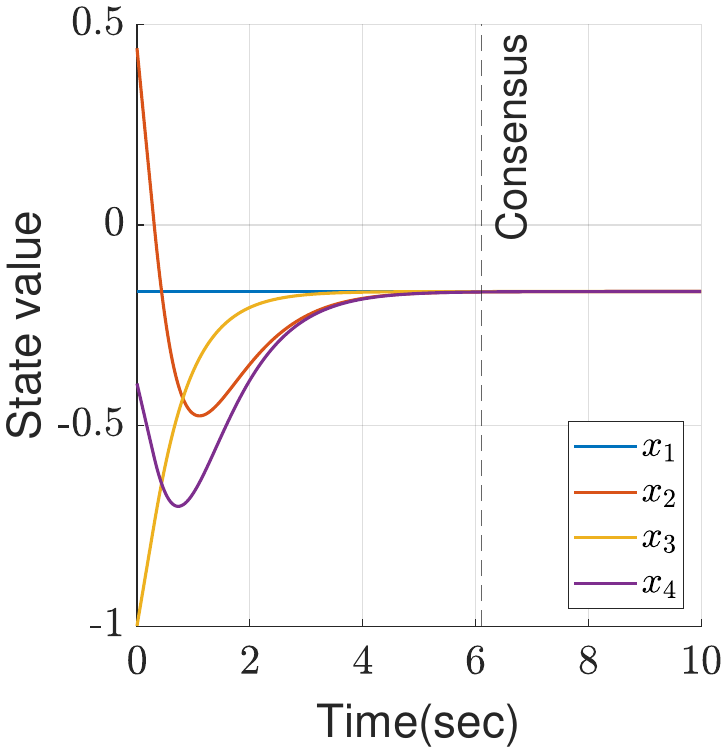}
		\caption{Consensus result of MAS in Fig. \ref{fig:mas}(\subref{fig:mas2}) with $\tau = 0.3s$.}
		\label{fig:consensus_result2}
	\end{subfigure}
	\hfill
	\begin{subfigure}{0.22\textwidth}
		\centering
		\includegraphics[width=0.99\textwidth]{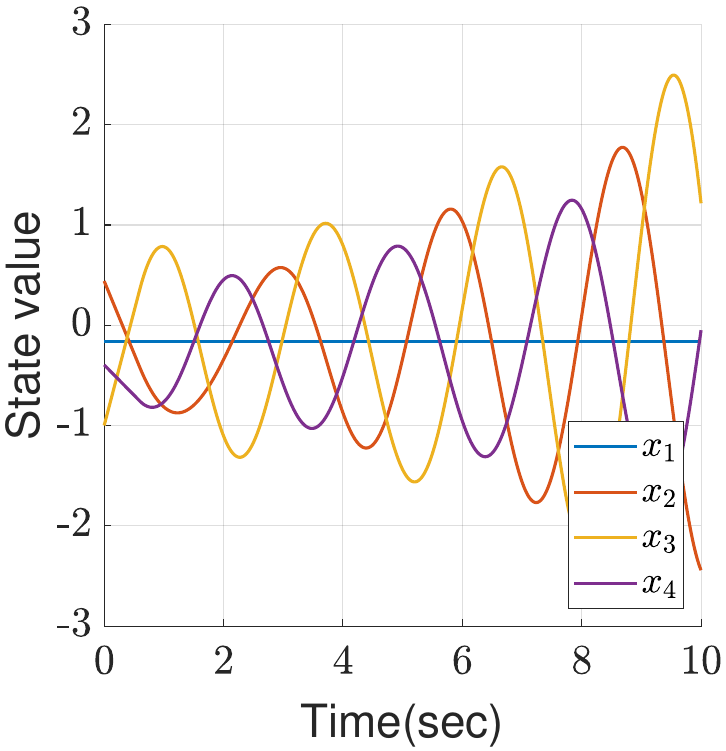}
		\caption{Consensus result of MAS in Fig. \ref{fig:mas}(\subref{fig:mas1}) with $\tau = 0.6s$.}
		\label{fig:consensus_result3}
	\end{subfigure}
	\hfill
	\begin{subfigure}{0.22\textwidth}
		\centering
		\includegraphics[width=0.99\textwidth]{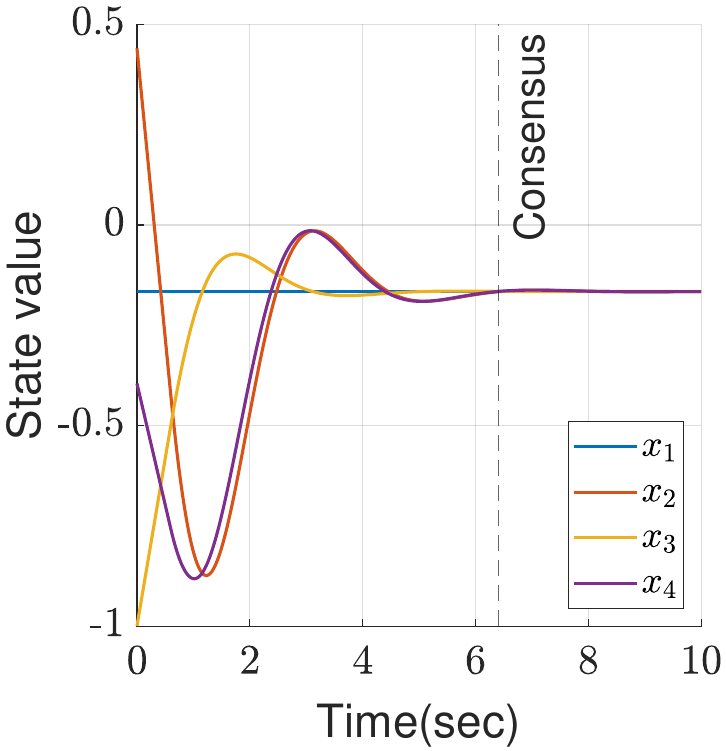}
		\caption{Consensus result of MAS in Fig. \ref{fig:mas}(\subref{fig:mas2}) with $\tau = 0.6s$.}
		\label{fig:consensus_result4}
	\end{subfigure}
	\caption{Consensus results of MASs.}
	\label{fig:consensus_result}
\end{figure}

Finally, let us provide a demo application of the obtained theoretical results via the consensus simulation of MASs. Consider a four-node MAS whose topology is modeled by an unweighted digraph shown in Fig. \ref{fig:mas}(\subref{fig:mas1}). The Laplacian eigenvalues of this digraph are (to two decimal places): $0$, $0.53$, $2.23 \pm 0.79\imath$, which contain complex values. The delayed consensus protocol of the MAS is as follows:
$$
\dot{x}(t)= -\Lg x(t - \tau)
$$
where $x(t) = [x_1(t), x_2(t), x_3(t), x_4(t)]^{\top}$, $x_1(t)$ to $x_4(t)$ are scalars describing the states of the four agents, $\tau \in \mathbb{R}^{+}$ is constant communication time-delay. Under the assumptions that the initial states of the agents are random values between $-1$ and $1$, $\tau = 0.3s$, the consensus result is illustrated in Fig. \ref{fig:consensus_result}(\subref{fig:consensus_result1}), where the consensus is numerically reached at $8.99s$ by taking $\max|x_i-x_j|<1e^{-3}$ as the threshold. However, if we remove the edge $(3,2)$ in Fig. \ref{fig:mas}(\subref{fig:mas1}) to obtain the digraph in Fig. \ref{fig:mas}(\subref{fig:mas2}) that satisfies the conditions of Theorem \ref{th:real_numbers} (new digraph has Laplacian eigenvalues: $0$, $1$, $1$, $2$), then the MAS numerically reaches consensus at $6.1s$ as shown in Fig. \ref{fig:consensus_result}(\subref{fig:consensus_result2}), with faster convergence and fewer oscillations compared to Fig. \ref{fig:consensus_result}(\subref{fig:consensus_result1}). Furthermore, if $\tau = 0.6s$, then the consensus protocol over the digraph in Fig. \ref{fig:mas}(\subref{fig:mas1}) will no longer be stable (see the trajectories in Fig. \ref{fig:consensus_result}(\subref{fig:consensus_result3})), while the consensus protocol over the digraph in Fig. \ref{fig:mas}(\subref{fig:mas2}) can still reach consensus at $6.41s$ (see the trajectories in  Fig. \ref{fig:consensus_result}(\subref{fig:consensus_result4})). It implies that the dynamical network of the digraph in Fig. \ref{fig:mas}(\subref{fig:mas2}) has higher tolerance to time-delay due to its purely real Laplacian spectrum.

\section{Conclusion}\label{sec:conclusion}

This paper has conducted a systematic investigation into the relationship between the topological properties and complex-valuedness of the Laplacian spectra of digraphs. It has been revealed that a real Laplacian spectrum is closely linked to the non-strong connectivity of a digraph, while the occurrence of complex Laplacian eigenvalues is related to the existence of directed cycles in a digraph. These findings are subsequently extended to multilayer digraphs.

This study opens a new vision for characterizing the Laplacian spectral properties of digraphs in terms of their topological properties, offering a pathway to ensure real Laplacian spectra of digraphs and mitigate the adverse effects of complex eigenvalues through topological design. Future works will explore the cost-effective ways to eliminate complex Laplacian eigenvalues and improve the performance of MASs via optimizing the topologies and edge weights of the underlying networks.

\ifCLASSOPTIONcaptionsoff
\newpage
\fi

{\footnotesize
	\bibliographystyle{IEEEtran}
	\bibliography{IEEEabrv,digraph-lap}
	
}

\end{document}